%% file: main.tex
\documentclass[twocolumn]{autart}
\usepackage{cite}
\usepackage{graphicx}
\usepackage{color}
\usepackage{textcomp}
\usepackage{amsmath}
\usepackage{booktabs}
\usepackage{amsfonts}
\usepackage{amssymb}
\usepackage{mathtools}
\usepackage{bm}
\usepackage{float}
\usepackage{comment}
\usepackage{subcaption}
\usepackage[normalem]{ulem}
\usepackage{multirow}
\usepackage[font={small}]{caption}
\usepackage{tabstackengine} 
\usepackage[mathscr]{euscript}
\usepackage{ulem}
\usepackage{bbm}
\graphicspath{ {images/} }
\usepackage[dvipsnames]{xcolor}

\usepackage{enumitem}

\begin{document}
\begin{frontmatter}
\input{sections/_title}
\author[UCSB-EE]{Shara Balakrishnan}\ead{sbalakrishnan@ucsb.edu},    
\author[UCSB-ME]{Aqib Hasnain}\ead{aqib@ucsb.edu},               
\author[PNNL]{Robert G.  Egbert}\ead{robert.egbert@pnnl.gov},  
\author[UCSB-ME]{Enoch Yeung}\ead{eyeung@ucsb.edu}              
\address[UCSB-EE]{Department of Electrical and Computer Engineering, University of California, Santa Barbara, United States} 
\address[UCSB-ME]{Department of Mechanical Engineering, University of California, Santa Barbara, United States} 
\address[PNNL]{Biological Sciences Division, Earth and Biological Sciences Directorate, Pacific Northwest National Laboratory, United States} 
\begin{abstract} 
\input{sections/0_abstract} 
\end{abstract}
\end{frontmatter}

\input{sections/1_Introduction}

\input{sections/2_BasicIdea}

\input{sections/3_MathematicalPreliminaries}
\input{sections/4_MainResults}

\input{sections/5_SimulatedResults}

\begingroup
\allowdisplaybreaks
\input{sections/6_Conclusion}
\input{sections/7_Acknowledgement}

\input{sections/8_Appendix}

\endgroup

\bibliographystyle{unsrt}
\bibliography{main,sensor_fusion,sensor_list,system_identification,koopman_operator_sensor_fusion_applications,biology_sensor_fusion,Koopman_delay_embedding}

\end{document}

%% file: sections/_title.tex
\title{Data-Driven Observability Decomposition with Koopman Operators for Optimization of Output Functions of Nonlinear Systems}

%% file: sections/0_abstract.tex
When complex systems with nonlinear dynamics achieve an output performance objective, only a fraction of the state dynamics significantly impacts that output. Those minimal state dynamics can be identified using the differential geometric approach to the observability of nonlinear systems, but the theory is limited to only analytical systems. In this paper, we extend the notion of nonlinear observable decomposition to the more general class of data-informed systems. We employ Koopman operator theory, which encapsulates nonlinear dynamics in linear models, allowing us to bridge the gap between linear and nonlinear observability notions. We propose a new algorithm to learn Koopman operator representations that capture the system dynamics while ensuring that the output performance measure is in the span of its observables. We show that a transformation of this linear, output-inclusive Koopman model renders a new minimum Koopman representation. This representation embodies only the observable portion of the nonlinear observable decomposition of the original system. A prime application of this theory is to identify genes in biological systems that correspond to specific phenotypes, the performance measure. We simulate two biological gene networks and demonstrate that the observability of Koopman operators can successfully identify genes that drive each phenotype. We anticipate our novel system identification tool will effectively discover reduced gene networks that drive complex behaviors in biological systems.



%% file: sections/1_Introduction.tex
\section{Introduction}
Sensor technology has advanced at a rapid pace, offering researchers unprecedented access to data on dynamical systems. Observability is the underlying principle that links the sensor data to the internal state of the system. Applications of observability include monitoring the state of the system \cite{rogne2018redundant,besanccon2007nonlinear,ortega2021generalized}, estimating process model parameters \cite{aguirre2018structural} and identifying optimal locations for sensor placement \cite{hinson2014observability}. The theory of observability is well established for linear systems \cite{hespanha2018linear}. Observability theory for nonlinear systems is limited to the differential geometric results for analytical systems\cite{nijmeijer1990nonlinear} and algebraic results for polynomial systems \cite{tibken2004observability}. For nonlinear systems learned from data, methods are being developed to identify if the system is observable\cite{yeung2018koopmanGramian}. The theory to identify the observable subspace decomposition of nonlinear systems from data-driven models is yet to be established.




Data-driven discovery of dynamics is critical for complex systems where the underlying mechanics are not fully understood. Such scenarios are common in biological cells \cite{gilpin2020learning}, finance \cite{chen2015modelling}, cyber-physical systems \cite{yuan2019data}, etc. One of the commonly used complex systems in biomanufacturing industries is the bacterium, \textit{Escherichia coli} \cite{kopp2020repetitive}. In   \textit{Escherichia coli}, gene transcription alone constitute over a $4,400-$dimensional dynamic process, and this excludes the protein and metabolic interactions within the cell. Such complex systems are typically deployed to achieve a specific performance objective. \textit{Escherichia coli} used in biomanufacturing processes are optimized for performance objectives like maximizing population cell growth \cite{shiloach2005growing} or maximizing production of a specified metabolite \cite{lin2006ethanol}. Only a fraction of the genes have a strong influence on the desired performance objective \cite{tong2020gene,wang2022screening,avwioroko2018isolation}.    This raises the question of how to identify a critical set of genes that have the strongest influence on given performance objective function. 

For a linear system, the performance objective can be treated as the output and an observable subspace decomposition results in the {\it minimal} system dynamics that drives the output \cite{zou2020moving}. Equivalent results have been developed for nonlinear systems using differential geometry for analytical systems where the governing equations are known prior \cite{nijmeijer1990nonlinear}. However, the dynamics of biological systems are not known prior and are typically learned from data. Hence, observable subspace decomposition methods cannot be used directly to learn the minimal gene expression dynamics in biological systems that drive a desired output phenotype.


In biological systems, the typical approach to identify genes that impact a phenotype is to look for genes that exhibit significant differences in their steady-state responses \cite{varet2016sartools,wetmore2015rapid,lima2013bacterial} across varying initial conditions. By considering initial conditions where the output (performance metric) response is vastly different, the genes with the highest differential steady state response are deemed to impact the output. This is a classical empirical approach that disregards both gene-to-gene interactions as well as gene-to-phenotype (output) interactions.  Our ultimate goal is to model these various nonlinear dynamical interactions from data and then find genes that drive a desired output which can later be used to optimize the performance of that output. 

Koopman operator theory is an increasingly popular approach to learn and analyze nonlinear system dynamics, specifically due to a growing suite of numerical methods that can be applied in a data-driven setting \cite{proctor2016dynamic,mezic2015applications}.  Koopman models are promising because they construct a set of state functions called Koopman observables that embed the nonlinear dynamics of a physical system in a high-dimensional space where the dynamics become linear \cite{budivsic2012applied}. Koopman models are typically learned from data using a dimensionality reduction algorithm called dynamic mode decomposition (DMD), which was developed by Schmid \cite{schmid2010dynamic}. Extensive research has enabled Koopman models to increase their predictive accuracy and decrease their computational complexity. Koopman models serve as a bridge between nonlinear systems and high-dimensional linear models, making them particularly helpful for extending linear notions to nonlinear systems in applications such as modal analysis\cite{mezic2013analysis,taira2017modal,mclean2020modal}, construction of observers \cite{surana2016linear,surana2016koopman,anantharaman2021koopman,yeung2018koopmanGramian,netto2018robust} and development of controllers\cite{proctor2016dynamic,korda2018linear,proctor2018generalizing,you2018deep,kaiser2021data}.

 The study of observability of nonlinear systems using Koopman operators is a growing area of research; Koopman operators have been augmented with output equations for applications like observer synthesis \cite{surana2016koopman,surana2016linear,anantharaman2021koopman}, optimal sensor placement\cite{hasnain2019optimalReporterPlacement,hasnain2022learning} and quantifying observability of nonlinear systems \cite{yeung2018koopmanGramian}. They all work under the assumption that the outputs lie in the span of Koopman observables but there is no theory on when that assumption holds. There are no algorithms to learn such output-inclusive Koopman models from data as Koopman models typically constitute a state equation learned either by using direct state measurements \cite{williams2015data,hasnain2019data,lusch2018deep} or delay-embedded output measurements \cite{bakarji2022discovering,balakrishnan2020prediction,arbabi2017ergodic}. Moreover, how to use Koopman operator models learnt from data to estimate the observable decomposition of the nonlinear system is yet to be established.
 

Here, we extend the theory of Koopman operators to nonlinear systems with a measurable output performance and develop the notion of observable subspaces for such nonlinear systems using linear Koopman operator theory. Through our investigation, we:
\begin{enumerate} [label=(\textit{\roman*})]
    \item developed a theory that maps the observable subspace of a nonlinear system to a linear output-inclusive Koopman model defined on that observable subspace (Theorems 3 and 4),
    \item identified the conditions under which the observable subspace of an output-inclusive Koopman model maps to the observable subspace of the nonlinear system (Theorem 5) 
    \item developed a new algorithm that learns such observable, output-inclusive Koopman models using deep learning and dynamic mode decomposition (Corollary 2), 
    \item showed that the new data-driven Koopman models can estimate the essential genes that drive the growth phenotype of a biological system in the order of their importance (Simulation Example 1), and
    \item showed that the gene dynamics in the observable subspace of each output of an interconnected genetic circuit constitute the significant genes that drive that output performance measure of the circuit (Simulation Example 2).
\end{enumerate}

The paper is organized as follows. Section \ref{sec: Basic Idea} introduces the problem statement in detail and Section \ref{sec: math prelims} briefly introduces the required mathematical preliminaries. In Section \ref{sec: Main Results}, we discuss the main theoretical results pertaining to observability of Koopman operators and the methods to see them in practice. We consider two simulated gene circuits in Section \ref{sec: Simulation Results} and demonstrate how the theory is used to find genes that drive each output of the system. Conclusions are drawn in Section \ref{sec: Conclusion}.

%% file: sections/2_BasicIdea.tex
\input{Images/FigureTexFiles/BasicIdea}

\section{Problem Formulation} \label{sec: Basic Idea}
We formulate the mathematical problem in more depth and describe how solving it benefits biological systems.
\subsection{The Mathematical Challenge} \label{subsec: Mathematical Formulation}
Given the autonomous discrete-time nonlinear dynamical system with output
\begin{subequations}\label{eq: nonlinear system with output}
    \begin{flalign}
        \text{State Equation:} && x_{t+1} &= f(x_t)&& \label{eq: nonlinear system without output}\\
        \text{Output Equation:} && y_t &= h(x_t)&& \label{eq: nonlinear output}
    \end{flalign}
\end{subequations}
where $x \in \mathcal{M} \subseteq \mathcal{R}^n$ is the state and $y \in \mathbb{R}$ is the output performance measure. The differential geometric approach to observability provides a nonlinear decomposition that can an analytical system of the form (\ref{eq: nonlinear system with output}) to 
\begin{align}\label{eq: nonlinear observable decomposition abstract}
    x^o_{t+1} &= f_o(x^o_t)\nonumber\\
    x^u_{t+1} &= f_u(x^o_t,x^u_t)\\
    y_t &= h_o(x^o_t)\nonumber
\end{align}
via a diffeomorphic (smooth and invertible) nonlinear transformation $\begin{bmatrix}x^o\\x^u\end{bmatrix} = \begin{bmatrix}\xi_o(x)\\\xi_u(x)\end{bmatrix} $ \vspace{3pt} where $x^u$ lies in the unobservable subspace of the system (\ref{eq: nonlinear system with output}). The remaining $x^o$ is the minimal state that drives the output dynamics and the manifold that $x^o$ lies in is the maximum subspace that the output $y$ can observe in the system (\ref{eq: nonlinear system without output}). We refer to that space observed by the output as the \textit{observable subspace} of the system (\ref{eq: nonlinear system with output}). For data-driven nonlinear models, there are no approaches to identify the nonlinear transformations $\xi_o$ and $\xi_u$. There are explicit methods to do similar transformations for data-driven linear systems and therefore, we turn to Koopman operator theory that  bridges the notions of linear and nonlinear observable decompositions.

A standard Koopman operator representation used to capture the nonlinear dynamical system with an output equation (\ref{eq: nonlinear system with output}) is given by 
\begin{subequations}\label{eq: Koopman operator system with output}
    \begin{flalign}
        \text{State Equation:} && \psi(x_{t+1}) &= K\psi(x_t)&& \label{eq: Koopman operator system dynamics}\\
        \text{Output Equation:} && y_t &= W_h\psi(x_t)&& \label{eq: Koopman operator output equation}
    \end{flalign}
\end{subequations}
where $\mathcal{M}\subseteq\mathbb{R}^n$ and $\psi:\mathcal{M}\rightarrow \mathbb{R}^{n_L}$ are the Koopman observables (functions of the state), whose linear evolution across time captures the nonlinear dynamics of the state and the output. To enable easier recovery of the base state $x$ from the Koopman observables $\psi(x)$, the Koopman observables are typically constrained to include the base states $x$ as
$\psi(x) = \begin{bmatrix} x^\top & \varphi^\top(x) \end{bmatrix}^\top$ where $\varphi(x)$ is a vector of pure nonlinear functions of $x$. The Koopman operators corresponding to the observables which contain the state $x$ are referred to as state-inclusive Koopman operators. For the rest of the paper, the Koopman model with observables denoted by $\psi$ are state-inclusive. Since the Koopman model (\ref{eq: Koopman operator system with output}) is linear, linear observability concepts can be used in this system. How do we use the Koopman system (\ref{eq: Koopman operator system with output}) to infer the observable state $x^o$ in (\ref{eq: nonlinear observable decomposition abstract})? Section \ref{sec: Main Results} delves more on this topic and provides algorithms to identify $x^o$ from data and determine the observable subspace of the original nonlinear system (\ref{eq: nonlinear system with output}).

\textbf{\subsection{The Biological Implication}}
In  complex microbial cell systems, techniques like transcriptomics \cite{wang2009rna} and proteomics \cite{aebersold2003mass} inform the dynamics within the cell and instruments like flow cytometers \cite{gant1993application}, plate readers \cite{meyers2018direct}, and microscopes \cite{lefman2004three} inform the phenotypic characteristics viewed from outside the cell. We can represent the intracellular activity by the state equation (\ref{eq: nonlinear system without output}) and the phenotype of interest by the output equation (\ref{eq: nonlinear output}). The phenotypic behavior is the performance metric that we wish to optimize with a specific objective. In Section \ref{sec: Simulation Results}, we simulate two biological gene networks, for which we learn the observable subspace of the nonlinear system (\ref{eq: nonlinear system with output}) and provide empirical methods to map that observation space (in which all of $x^o$ lies) to the set of genes (a subset of the state variables in $x$) that drive the output phenotypic behavior. Upon identifying the genes that influence the phenotypic dynamics, we can deploy actuators developed for biological systems to control the gene expression and optimize the phenotypic performance. 

The generic phenotypic performance optimization problem is given by 
\begin{align}
    \hspace{-1in}
    &\quad \max_{u} \sum_{t=0}^N||y_t||^2_2\\
    \text{such that }\quad \notag   & x_{t+1} = f(x_t) + \sum_{i=1}^{n_a} g_i(x_t,u_{t,i})\notag \\
    & y_t = h(x_t) \notag
\end{align}
where $g$ is the input function that captures both how an input directly controls the expression of targeted genes as well as off-target gene expression effects \cite{chen2018model} and $n_a$ is the number of individual genes whose expression dynamics we can target to control. Two accessible genetic actuators that control gene expression are: A) Transposons \cite{wetmore2015rapid} which knockout the complete gene expression with $g_i(x_t,u_{t,i}) = -f_i(x_t)$, and B) CRISPR interference mechanism which suppresses the gene expression \cite{huang2021dcas9} with $g_i(x_t,u_{t,i}) <0$. We anticipate this work will enable the identification of genes that impact growth of soil bacteria in sparse environmental conditions that can be controlled  by biological actuators to maximize their population growth.

%% file: Images/FigureTexFiles/BasicIdea.tex
\begin{figure*}[htp!]
\centering
    \includegraphics[width=\textwidth*2/3,keepaspectratio]{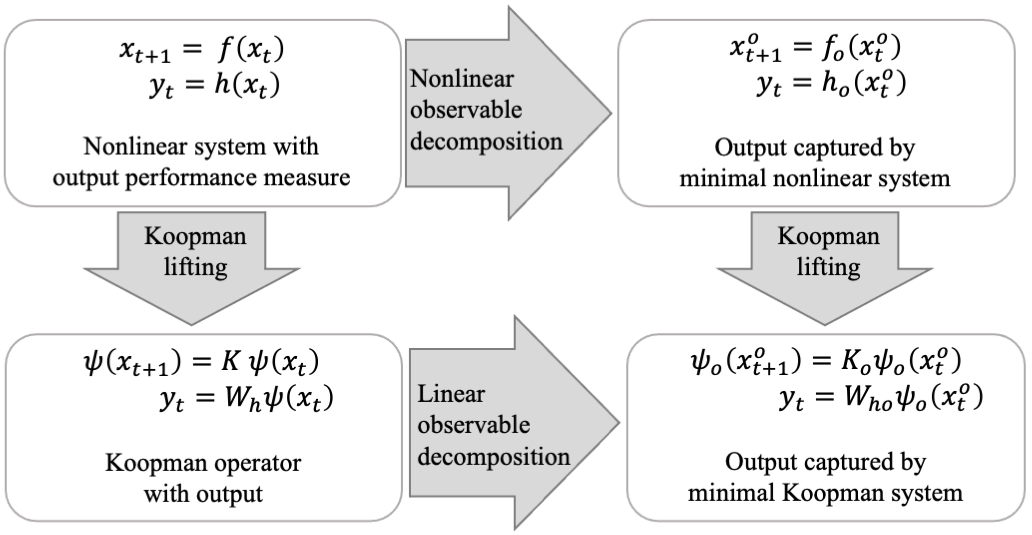}
    \caption{\textbf{Koopman approach to observability decomposition of nonlinear systems:} The nonlinear observable decomposition (upper transition) is a result from the the differential geometric approach to observability of nonlinear systems which is only defined for analytical systems. The Koopman lifting (transition on the left) is from Koopman operator theory to find high-dimensional linear representations of nonlinear system. Our approach is to find the structure of the Koopman operator for the nonlinear decomposed system (transition on the right) and establish a relationship with the Koopman operator model of the original nonlinear system through a linear transformation (lower transition).}\label{fig: example1}
\end{figure*}

%% file: sections/3_MathematicalPreliminaries.tex
\section{Mathematical Preliminaries}\label{sec: math prelims}
We consider the discrete-time nonlinear dynamical system (\ref{eq: nonlinear system with output}) with  the state $x \in \mathcal{M} \subseteq \mathbb{R}^n$ and an output $y \in \mathbb{R}$. The Koopman operator for the state dynamics (\ref{eq: nonlinear system without output}) is given by (\ref{eq: Koopman operator system dynamics}) where $K:\mathcal{F}^{n_L} \rightarrow \mathcal{F}^{n_L} $, $\psi:\mathcal{M} \rightarrow \mathcal{R}^{n_L}$, and $ \mathcal{F}$ is a space of smooth functions. The Koopman operator is ideally a linear infinite dimensional operator($n_L \rightarrow \infty$) but Koopman models identified from data are finite dimensional approximations ($n_L<\infty$). Detailed discussion on Koopman operator theory can be found in \cite{mauroy2020koopman,brunton2021modern}. Since the core focus of the paper is on the observability of (\ref{eq: nonlinear system with output}), we present the relevant results from the differential geometric approach to observability of discrete-time nonlinear systems \cite{nijmeijer1982observability, nijmeijer1990nonlinear,hanba2009uniform,albertini1996remarks}. In addition, we also present an overview of the existing algorithms used to identify Koopman operators.

\subsection{Observability of discrete-time nonlinear systems} \label{subsec: NL observability}
The observability of the nonlinear system (\ref{eq: nonlinear system with output}) revolves around the properties of a new space obtained by the transformation of the base coordinates $x$, called the \textit{observation space}. 

\begin{defn}
The observation space $\mathcal{O}_y(x)$ for the nonlinear dynamical system (\ref{eq: nonlinear system with output}) is the space of functions that captures the output across infinite time:
\begin{align*}
    \mathcal{O}_y(x) = \{h(x), h(f(x)), \cdots, h(f^i(x)), \cdots\}, \quad i \in \mathbb{Z}_{>0}.
\end{align*}
\end{defn}

With a slight abuse of notation, based on the context, we use $\mathcal{O}_y(x)$ to represent either a set or a vector of functions. If the observation space $\mathcal{O}_y(x)$ has a diffeomorphic map (smooth and invertible) with $x$, then the outputs across infinite time can be used to estimate the initial state $x$ and this would be true for all $x \in \mathcal{M}$. This is the strongest condition that ensures the system (\ref{eq: nonlinear system with output}) is observable, but it is impossible to check for. So, a more local approach is adopted by computing the dimension of the observation space at a point. 

\begin{defn}
The dimension of the observation space $\mathcal{O}_y(x)$ at a point $\bar{x} \in \mathcal{M}$ is the rank of the Jacobian matrix of the function set $ \{h(x), h(f(x)), \cdots, h(f^{n-1}(x))\}$:
\[
dim\Big{(}\mathcal{O}_y(\bar{x})\Big{)} \hspace{-3pt}=\hspace{-1pt}rank \hspace{-3pt}
\left.
\begin{bmatrix}
    \frac{\partial h(x)}{\partial x_1} & \cdots & \frac{\partial h(x)}{\partial x_n} \\
    \vdots  & \ddots & \vdots\\
    \frac{\partial h(f^{n-1}(x))}{\partial x_1}  & \cdots & \frac{\partial h(f^{n-1}(x))}{\partial x_n}
\end{bmatrix}\right\vert_{x=\bar{x}}.
\]
\end{defn}

The dimension of the observation space can be computed locally at a point and hence a local observation result can be obtained. While there are different notions of observability for nonlinear systems, we only discuss strongly local observability as we build on top of this definition for the rest of the paper.

\begin{defn}
The system (\ref{eq: nonlinear system with output}) is said be strongly locally observable at $x\in \mathcal{M}$ if there exists a neighborhood $\mathcal{U}$ of $x$ such that for any $\bar{x} \in \mathcal{U}$,  $h(f^{k}(\bar{x})) = h(f^{k}(x))$ for $k=0,1,\cdots,n-1$, implies $\bar{x} =x$. 
\end{defn}

\begin{thm}\label{theorem: nonlinear observability}
\big{(}Theorem 2.1 from \cite{nijmeijer1982observability}\big{)}
If the system (\ref{eq: nonlinear system with output}) is such that $dim(\mathcal{O}_y(\bar{x}))=n$,  then the system is strongly locally observable at $\bar{x}$
\end{thm}
The results extend to the full system if they are true for all $x\in \mathcal{M}$. In that case, we state that the system is strongly locally observable if $dim(\mathcal{O}_y(x))=n$ $\forall x \in \mathcal{M}$. The premise of the paper is that for complex dynamics, the system is not strongly locally observable and that only a subspace of the system is observable. This is captured in part as continuous and discrete results in:
Remark 2 under Theorem 2.9 under \cite{nijmeijer1982observability}, Proposition 3.34 from \cite{nijmeijer1990nonlinear} and Theorem 3.51 from \cite{nijmeijer1990nonlinear}. We consider this result for discrete systems in a form that is useful to us.

\begin{thm} \label{theorem: prelim nonlinear observable decomposition}
Given a nonlinear system (\ref{eq: nonlinear system with output}), for a given point $x \in \mathcal{M}$, if there exists a neighborhood $\mathcal{U}$ of $x$ such that for any $\bar{x} \in \mathcal{U}$, $dim(\mathcal{O}_y(\bar{x}))=r$, then we can find a local coordinate transform from $x$ to $\tilde{x} = \begin{bmatrix}
    \tilde{x}^1 \\ \tilde{x}^2
\end{bmatrix}$ such that
\begin{subequations}\label{eq: nonlinear observable decomposition}
    \begin{align}
        \tilde{x}^1_{t+1} &= f_o(\tilde{x}^1_t) \label{eq: nonlinear observable decomposition state}\\
    \tilde{x}^2_{t+1} &= f_u(\tilde{x}^1_t,\tilde{x}^2_t)\\
    y_{t} &= h_o(\tilde{x}^1_t)\label{eq: nonlinear observable decomposition output}
    \end{align}
\end{subequations}
 $\forall$ $x \in \mathcal{U}$ where $\tilde{x}^1 \in \mathcal{M}'\subseteq \mathbb{R}^r$ and $\mathcal{M}'\subseteq \mathcal{M}$.
\end{thm}
The result tells us that if the output dynamics (given by $\mathcal{O}_y(x)$) lies in a lower dimensional space than the state dynamics, then we can do a nonlinear observable decomposition of the original system(\ref{eq: nonlinear system with output}) to get a system of the form (\ref{eq: nonlinear observable decomposition abstract}) where only a subset of the states are observable and drive the output dynamics. This is an important result as it is this system that needs to be connected to the Koopman model (\ref{eq: Koopman operator system with output}) to arrive at the observability analysis of the nonlinear system (\ref{eq: nonlinear system with output}) when identified from data.

\subsection{Dynamic mode decomposition (DMD)} \label{subsec: DMD}
Dynamic Mode Decomposition (DMD) is a popular class of algorithms adopted to learn approximate finite dimensional Koopman models. A detailed review of some of the popular DMD algorithms is given in \cite{schmid2021dynamic}. A general framework for the DMD algorithm is given by
\begin{align} \label{eq: general DMD formulation}
\min_{\psi,K}&||\psi(x_{k+1})-K\psi(x_k)||^2_F
\end{align}
where the number of Koopman observables ($n_L\geq n$) is a typical hyperparameter. The exact DMD algorithm \cite{schmid2010dynamic} identifies local linear representations of (\ref{eq: nonlinear system without output}) by setting $\psi(x) =x$. The extended DMD (E-DMD) algorithm proposed in \cite{williams2015data} identifies Koopman models by using a kernel of user-specified functions to represent the Koopman observables $\psi(x)$. To automate the Koopman observable learning process, deepDMD algorithms \cite{yeung2019learning,li2017extended} specify $\psi(x)$ as the output layer of a neural network:
\begin{align*} \label{eq: deep DMD formulation}
\psi(x) = 
\begin{bmatrix}
    x\\
    \varphi(x)
\end{bmatrix} =& 
\begin{bmatrix}
    x\\
    g_n \circ \sigma \circ \cdots \circ \sigma \circ g_2 \circ \sigma \circ g_1(x)
\end{bmatrix}\nonumber
\end{align*}
where the $i^{th}$ hidden layer captured by weights $W_i$, biases $b_i$, linear function $g_i(x) = W_ix+b_i$ and activation function $\sigma$. Activation functions like sigmoidal \cite{cybenko1989approximation}, rectified linear unit (ReLU) activation functions  \cite{hanin2019universal} and radial basis functions (RBFs) \cite{lo1998multilayer} parameterize $\psi(x)$ with universal function approximation properties. There are other algorithms which identify non state-inclusive Koopman operators by identifying purely nonlinear Koopman observables that have a differomorphic map with the base state $x$. The diffeomorphic map is implemented using autoencoders \cite{otto2019linearly,takeishi2017learning}. 


%% file: sections/4_MainResults.tex
\section{Main Results}\label{sec: Main Results}
In this section, we methodically show how we translate the theory of observable decomposition of nonlinear systems to the linear observable decomposition of Koopman systems to discover the critical Koopman observables(functions of the state $x$) that drive the output dynamics.  We use an analytical example to elucidate our theoretical results. Along the way, we discuss how to extend the theory to practise. The details for implementing the algorithms are provided in the Appendix (Section \ref{sec: appendix}). 

\subsection{Minimal Koopman operator that drives the output}
We begin by showing the nonlinear system (\ref{eq: nonlinear system with output}) can be transformed into a minimal Koopman model that drives the output performance metric.
\input{Theorems/Theorem_1_r-dimensional-reduction}
\input{Theorems/Proofs/Proof_Theorem1}

The nonlinear observable decomposition theorem transforms the full nonlinear system to a new coordinate space with the minimal number of state variables  required to capture the output performance metric. Theorem \ref{Theorem: Theoretical r-dimensional OC-KO existence} proves the existence of a state-inclusive Koopman operator for the transformed system such that its Koopman observables lie in the span of the observation space vector and vice versa. The above theorem is fusing the information in the states and outputs. One of the existing results in fusing two measurements is the theory of factor conjugacy in \cite{williams2015datafusion,mezic2019spectrum}. In the following remark, we tie Theorem 1 to the concept of factor conjugacy.

\begin{rem} 
If there are two dynamical systems $x_{t+1} = f_x(x_t)$ and $z_{t+1} = f_z(z_t)$ such that $z = h_{zx}(x)$, then the two systems are said to be factor conjugate if $h_{zx}(f_x(x)) = f_z(h_{zx}(x))$. As a consequence, if $\{(\lambda_{zi},\phi_{zi}), i=1,2,...\}$ represents the set of all eigenvalue-eigenfunction pair of the Koopman operator for the $z-$dynamics, then $\{(\lambda_{zi},\phi_{zi}\circ h_{zx}), i=1,2,...\}$ represent a \textit{subset} of all eigenvalue-eigenfunction pair of the Koopman operator for the $x-$dynamics. The subset of eigenfunctions $\{\phi_{zi}\circ h_{zx}, i\in \mathbb{N}\}$ in the $x$-dynamics are a minimal set of observables required to capture the output dynamics and hence constitute a basis for the reuduced Koopman operator (\ref{eq: Koopman operator of Sussman decomposition}).
\end{rem}

In certain nonlinear systems, all functions in the observation space $\mathcal{O}_y(x)$ lie in the span of a finite subset of functions in  $\mathcal{O}_y(x)$. This strong criteria results in finite dimensional exact Koopman operators. This is a useful result for elucidating Theorem \ref{Theorem: Theoretical r-dimensional OC-KO existence} later and is formally stated in the below corollary.

\begin{cor}
If there exists a finite dimensional observation space $\mathcal{O}_{y,q}(x) = \{h(x), h(f(x)), \cdots, h(f^q(x))\}$ such that any function $h(f^i(x)) \in \mathcal{O}_y(x)$ where $i \in \mathbb{Z}_{\geq0}$ lies in the span of $\mathcal{O}_{y,q}(x)$ and dim$(\mathcal{O}_{y,q}(x))=r \leq n$, we can find a finite dimensional exact Koopman operator representation of the form (\ref{eq: Koopman operator of Sussman decomposition}). 
\end{cor}

\subsection{Learning Koopman operators with output}\label{sec: OC-DMD}
We explored how dynamic mode decomposition (DMD) algorithms are used to learn approximate Koopman operators in Section \ref{subsec: DMD}. In prior works where an output equation is involved with Koopman operators \cite{mesbahi2021nonlinear,yeung2018koopmanGramian,surana2016linear,surana2016koopman,netto2018robust}, the outputs are typically assumed to lie in the span of the Koopman observables; there are no DMD algorithms to ensure that. The following corollary to Theorem \ref{Theorem: Theoretical r-dimensional OC-KO existence} relaxes that assumption and provides the necessary and sufficient condition for the existence of a Koopman operator representation of form (\ref{eq: Koopman operator system with output}).
\begin{cor} \label{cor: necessary and sufficient condition for OC-DMD}
Given the dynamical system (\ref{eq: nonlinear system with output}), we can find the output-constrained Koopman operator representation (\ref{eq: Koopman operator system with output}) if and only if the observation space $\mathcal{O}_y(x)$ of (\ref{eq: nonlinear system with output}) lies in the span of the Koopman observables. 
\end{cor}

  We incorporate Corollary \ref{cor: necessary and sufficient condition for OC-DMD} into the DMD objective \ref{eq: general DMD formulation} to form the more generic DMD multi-objective optimization problem:
\begin{align*}
    \min_{\psi,K,W_h} ||\psi(x_{t+1}) &- K\psi(x_t)||_F^2 + ||y_t-W_h\psi(x_t)||_F^2
\end{align*}
where the output is forced to lie in the span of the observables. Moreover, if the output at time $t$ lies in the span of the observables, i.e., $y_t = W_h\psi(x_t)$, then for any future time point $t+k$, the output at that time point also lies in the span of the observables $\psi(x_t)$ as $y_{t+k} = W_hK^k\psi(x_t)$. This ensures that the full observation space $\mathcal{O}_y(x)$ of the nonlinear system (\ref{eq: nonlinear system with output}) lies in the span of the Koopman observables $\psi(x)$, thereby adhering to Corollary 2. Since the above objective function forces the output to lie in the span of the Koopman observables, we term this problem as \textit{Output constrained dynamic mode decomposition} (OC-DMD). The neural network based implementation of OC-DMD is termed as OC-deepDMD, the details of which is discussed in Appendix \ref{appendix: ocdeepDMD}.

\subsection{Identifying the minimal Koopman operator}
Theorem \ref{Theorem: Theoretical r-dimensional OC-KO existence} establishes the existence of a minimal Koopman operator model that drives the output performance metric. How do we learn this model in practice? The following result establishes a procedure to do so.

\input{Theorems/Theorem_2_OCKOR_to_NLObservableDecomposition}
\input{Theorems/Proofs/Proof_Theorem2}

Theorem \ref{Theorem: Linear Observable Decomposition of OCKO} provides a route to identify the minimal Koopman operator (\ref{eq: Koopman operator of Sussman decomposition}) that drives the output performance metric of the nonlinear system (\ref{eq: nonlinear system with output}). We can use the OC-DMD algorithm from Section \ref{sec: OC-DMD} to identify a Koopman operator representation with an output equation (\ref{eq: Koopman operator system with output}) and then use a linear transformation $\psi_o(x^o)=T\psi(x)$ to go from (\ref{eq: Koopman operator system with output}) to (\ref{eq: Koopman operator of Sussman decomposition}). The next obvious question is--- what is the linear transformation $T$? We use the observable decomposition approach in Linear systems \cite{hespanha2018linear} to find $T$ for nonlinear systems with analytical finite dimensional Koopman operator representations.

\begin{cor} \label{cor: linear observable decomposition}
Suppose the nonlinear system (\ref{eq: nonlinear system with output}) has an exact finite dimensional Koopman operator representation (\ref{eq: Koopman operator system with output}) and a minimal finite dimensional Koopman operator representation of the form (\ref{eq: Koopman operator of Sussman decomposition}) where $x_o =\xi_o(x)$,  $x_o \in \mathcal{M}'\subseteq \mathbb{R}^r$ and $\psi_o(x^o):\mathcal{M}'\rightarrow \mathbb{R}^{{n}_{oL}}$. If $V$ represents the matrix of right singular vectors of the observability matrix of (\ref{eq: Koopman operator system with output})
$
    \mathcal{O}_\psi = 
    \begin{bmatrix}
        W_h^\top & K^\top W_h^\top & \cdots & (K^{n_L})^\top W_h^\top 
    \end{bmatrix}^\top,
$
then the transformation $\begin{bmatrix}
    \psi_1(x)\\
    \psi_2(x)
\end{bmatrix} = V^\top\psi(x)$ results in the observable decomposition form
\begin{align*}
    \begin{bmatrix}
        \psi_1(x_{t+1})\\
        \psi_2(x_{t+1})
    \end{bmatrix}
    &=
    \begin{bmatrix}
        K_1 & 0\\
        K_{12} & K_2
    \end{bmatrix}
    \begin{bmatrix}
        \psi_1(x_{t})\\
        \psi_2(x_{t})
    \end{bmatrix}\\
    y_t &= \begin{bmatrix}
        W_{h1} & 0
    \end{bmatrix}
    \begin{bmatrix}
        \psi_1(x_{t})\\
        \psi_2(x_{t})
    \end{bmatrix}
\end{align*}
where $\begin{bmatrix}
        K_1 & 0\\
        K_{12} & K_2
    \end{bmatrix} = V^\top K V$, $\begin{bmatrix}
        W_{h1} & 0
    \end{bmatrix} = W_hV$ and $\psi_1(x)$ are the state variables $\psi_o(x^o)$ of the minimal Koopman operator (\ref{eq: Koopman operator of Sussman decomposition}), i.e., $\psi_1(x)=\psi_o(x^o)$. 
\end{cor}

The above corollary provides a method to identify the minimal Koopman operator representation that drives the output dynamics. This is the same approach that we can adopt in practice for the approximate finite dimensional Koopman operators learnt from OC-DMD algorithm. The details of this approach is discussed in Appendix \ref{appendix: observable decomposition}. The uniqueness of the solution is discussed in the following remark.

\begin{rem}\label{remark: x_o not unique }
Corollary 2 provides a way to transform (\ref{eq: Koopman operator system with output}) to (\ref{eq: Koopman operator of Sussman decomposition}) but it does not provide the exact expression of $x^o$. The reason is that $x^o$ is not unique; any set of $r$ functions selected from $\psi_1(x)$ having a Jacobian of rank $r$ is a valid candidate for $x^o$. All functions of $\psi_1(x)$ can be written as either a linear or a nonlinear function of that $x^o$.
\end{rem}

\subsection{State information contained in the outputs}

An important practical consideration in complex systems is to gauge if the sensor measurements obtained from the system (\ref{eq: nonlinear system with output}) constitute a representation of the system state; in other words, does the fusion of the sensor measurements have a diffeomorphic relationship with the system state $x$. We use the established concepts in observability analysis to answer that question in the following theorem.

\input{Theorems/Theorem_3_Fusion_of_multiple_outputs}
\input{Theorems/Proofs/Proof_Theorem3}

Theorem \ref{Theorem: delay embedded Koopman operator} provides a framework to check if the fusion of various output measurements across space and time renders a representation of a state for the complex system dynamics (\ref{eq: nonlinear system with output}). The observability angle provides insight on why delay embedded Koopman observables are useful in the identification of Koopman operators \cite{arbabi2017ergodic,bakarji2022discovering, balakrishnan2020prediction}.

\begin{rem}
The Koopman observables of the delay embedded Koopman operator in Theorem \ref{Theorem: delay embedded Koopman operator}, $\psi(z)$ is the union of the observation spaces of all the individual output measurements \big{(}$\mathcal{O}_{y_1}(x) \cup \mathcal{O}_{y_2}(x) \cup \cdots \cup \mathcal{O}_{y_p}(x) $\big{)}. 
\end{rem}

In the observable subspace identification problem, one of the concerns to be wary of is whether or not the outputs have sufficient information about that subspace. Theorem \ref{Theorem: delay embedded Koopman operator} can be used to learn the delay embedded Koopman operator to capture the dynamics of (\ref{eq: nonlinear system with output}) and check if there is a diffeomorphism between $\psi(z)$ and $x$. The detailed procedure is given in Appendix \ref{appendix: diffeomorphic map}. The reason why the full $\psi(z)$ is used and not just $z$ is that $x^o$ is not unique (as seen in Remark \ref{remark: x_o not unique }) and any $n$ functions in $\psi(z)$ could form a diffeomorphic map with $x$.

\input{Examples/Theoretical_Example_1}

%% file: Theorems/Theorem_1_r-dimensional-reduction.tex
\begin{thm}\label{Theorem: Theoretical r-dimensional OC-KO existence}
Suppose the dynamical system with the output performance measure (\ref{eq: nonlinear system with output})
\begin{align*}
    x_{t+1} &= f(x_t)\\
    y_t &=h(x_t)
\end{align*}
where $x \in \mathcal{M} \subseteq \mathbb{R}^n$ and $y \in \mathbb{R}$ is such that its observation space $\mathcal{O}_y(x)$ has a constant dimension $r \leq n$ at $x$ in the neighborhood $\mathcal{U} \subseteq \mathcal{M}$. Then, for any $x \in \mathcal{U}$ there are $r$ functions in $\mathcal{O}_y(x)$ which constitute a surjective coordinate transformation to a reduced space $x^o = \xi(x) \in \mathcal{M}' \subseteq \mathbb{R}^r$  with a Koopman operator representation
\begin{equation}\label{eq: Koopman operator of Sussman decomposition}
    \begin{aligned}
        \psi_o(x^o_{t+1}) &= K_o \psi_o(x^o_{t})\\
        y_t & = W_{ho}\psi_o(x^o_t).
    \end{aligned}
\end{equation}
\end{thm}

%% file: Theorems/Proofs/Proof_Theorem1.tex
\begin{pf} The proof involves three steps:
\vspace{0pt}

(\textit{i}) \textit{Convert the system to the nonlinear observable canonical form:}
It is given that at a point $x$ in a neighborhood $\mathcal{U}$, we have $dim(\mathcal{O}_y(x)) = r\leq n$. Using Theorem \ref{theorem: prelim nonlinear observable decomposition}, we can transform the the base state $x$ to a new local coordinate $(\tilde{x}^1,\tilde{x}^2)$ \vspace{3pt} using a diffeomorphic local coordinate transform 
$
    \begin{bmatrix}
    \tilde{x}^1\\
    \tilde{x}^2
    \end{bmatrix} = 
    \begin{bmatrix}
    \xi_o(x)\\
    \xi_u(x)
    \end{bmatrix}
$\vspace{3pt}
such that (\ref{eq: nonlinear system with output}) can take the nonlinear observable canonical form (\ref{eq: nonlinear observable decomposition}). 
\vspace{0pt}

(\textit{ii}) \textit{Find an infinite dimensional linear model to represent the nonlinear observable canonical form:} Consider the observation space of (\ref{eq: nonlinear observable decomposition}) in the vector form:
\[
\bar{\mathcal{O}}_y(\tilde{x}^1) \triangleq \begin{bmatrix}h_o(\tilde{x}^1)^\top & h_o(f_o(\tilde{x}^1))^\top & h_o(f_o^2(\tilde{x}^1))^\top & \cdots \end{bmatrix}^\top.
\]
When we propagate $\bar{\mathcal{O}}_y(\tilde{x}^1)$ from time point $t$ to $t+1$, all functions of $\bar{\mathcal{O}}_y(\tilde{x}^1_{t+1})$ lie in the span of all functions in $\bar{\mathcal{O}}_y(\tilde{x}^1_{t})$ as $h_o(f_o^i(\tilde{x}^1_{t+1})) = h_o(f_o^{i+1}(\tilde{x}^1_{t}))$ $\forall$ $i \in \mathbb{Z}_{\geq 0}$. Hence, there exists an infinite dimensional matrix $K_y$ that renders an infinite dimensional linear model:
\begin{align*}
    \bar{\mathcal{O}}_y(\tilde{x}^1_{t+1}) &= K_y \bar{\mathcal{O}}_y(\tilde{x}^1_{t})\\
    y_t &= \begin{bmatrix}\mathbb{I}_p & 0\end{bmatrix}\bar{\mathcal{O}}_y(\tilde{x}^1_{t})
\end{align*}
for the nonlinear dynamics of $\tilde{x}^1$ given by (\ref{eq: nonlinear observable decomposition state}) and (\ref{eq: nonlinear observable decomposition output}). 
\vspace{0pt}

(\textit{iii}) \textit{Convert the infinite-dimensional system to a state-inclusive Koopman operator representation:} For output at a time point $t+k$, we have from systems (\ref{eq: nonlinear system with output}) and (\ref{eq: nonlinear observable decomposition}), $y_{t+k} = h(f^k(x_t)) = h_o(f_o^k(\tilde{x}^1_t))$ for any $x_t \in \mathcal{M}$. Hence, we have $\mathcal{O}_y(x) =\bar{\mathcal{O}}_y(\tilde{x}^1) \Rightarrow dim(\bar{\mathcal{O}}_y(\tilde{x}^1)) =r$. So, using Theorem 2.8 in \cite{nijmeijer1982observability} or the discrete time equivalent of Proposition 3.34 in \cite{nijmeijer1990nonlinear} , we can find a vector of $r$ functions in the observation space $\bar{\mathcal{O}}_y(\tilde{x}^1)$, say $x^o\in \mathbb{R}^r$, such that all functions in $\bar{\mathcal{O}}_y(\tilde{x}^1)$ are a nonlinear function of $x^o$. Then, there exists a permutation matrix $P$ such that $P\bar{\mathcal{O}}_y(\tilde{x}^1) = \begin{bmatrix}
    x^o\\
    \varphi_o(x^o)
\end{bmatrix} = \psi_o(x^o)$\vspace{2pt} which converts the above infinite-dimensional linear model to a state-inclusive Koopman operator representation of the form (\ref{eq: Koopman operator of Sussman decomposition}). Hence the proof.
\hspace*{\fill}~\qed\par
\end{pf}

%% file: Theorems/Theorem_2_OCKOR_to_NLObservableDecomposition.tex
\begin{thm} \label{Theorem: Linear Observable Decomposition of OCKO}
Suppose the dynamical system with the output performance measure (\ref{eq: nonlinear system with output})
\begin{align*}
    x_{t+1} &= f(x_t)\\
    y_t &=h(x_t)
\end{align*}
where $x \in \mathcal{M} \subseteq \mathbb{R}^n$ and $y \in \mathbb{R}$ is such that its observation space $\mathcal{O}_y(x)$ has a constant dimension $r \leq n$ at $x$ in the neighborhood $\mathcal{U} \subseteq \mathcal{M}$. For $x \in \mathcal{U}$, if the nonlinear system (\ref{eq: nonlinear system with output}) 
has a Koopman operator representation (\ref{eq: Koopman operator system with output})
\begin{align*}
    \psi(x_{t+1}) &= K\psi(x_t)\\
    y_t &= W_h\psi(x_t),
\end{align*}
then there exists a linear coordinate transform $T$ that takes the Koopman operator (\ref{eq: Koopman operator system with output}) to the minimal Koopman operator (\ref{eq: Koopman operator of Sussman decomposition}) that drives the output performance. 
\end{thm}

%% file: Theorems/Proofs/Proof_Theorem2.tex
\begin{pf}
For $x \in \mathcal{U}$, it is given that the nonlinear system (\ref{eq: nonlinear system with output}) has a Koopman operator representation of the form (\ref{eq: Koopman operator system with output}). By corollary \ref{cor: necessary and sufficient condition for OC-DMD}, it is evident that the observation space $\mathcal{O}_y(x)$ of the nonlinear system (\ref{eq: nonlinear system with output}) lies in the span of the Koopman observables, i.e., there exists a transformation $T_1$ such that $\mathcal{O}_y(x) = T_1\psi(x)$. Since  $dim(\mathcal{O}_y(x)) =r$, using \ref{Theorem: Theoretical r-dimensional OC-KO existence} we can get the minimal Koopman operator model (\ref{eq: Koopman operator of Sussman decomposition}) with the property that the Koopman observables $\psi_o(x^o)$ lie in the span of the observation space $\mathcal{O}_y(x)$ ($=\bar{\mathcal{O}}_y(\tilde{x}^1)$). Hence, there exists a linear transformation $T_2$ such that $\psi_o(x^o) =T_2\mathcal{O}_y(x)$. Therefore, there exists a linear coordinate transform $T=T_2T_1$ that takes the full Koopman operator representation (\ref{eq: Koopman operator system with output}) to the minimal Koopman operator representation (\ref{eq: Koopman operator of Sussman decomposition}) that drives the output performance metric. Hence the proof. 
\hspace*{\fill}~\qed\par
\end{pf}


%% file: Theorems/Theorem_3_Fusion_of_multiple_outputs.tex
\begin{thm} \label{Theorem: delay embedded Koopman operator}
Given that the $n-$dimensional nonlinear dynamical system (\ref{eq: nonlinear system with output}) with $p$ output measurements:
\begin{align*}
    x_{t+1} &= f(x_t)\\
    y_t &= h(x_t)
\end{align*}
has an observation space $\mathcal{O}_y(x)$ of constant dimension $r$ for all points $x\in \mathcal{U}$ where $\mathcal{U} \subseteq \mathcal{M}$. Then, there exists $n_d \in \mathbb{N}$ such that the delay embedded output
\[
z_{t} = \begin{bmatrix}
    y_{n_dt}^\top & y_{n_dt+1}^\top \cdots & y_{n_d(t+1)-1}^\top
\end{bmatrix}^\top
\]
has a Koopman operator representation
\begin{align}\label{eq: delay embedded Koopman operator}
    \psi_z(z_{t+1}) &= K_z\psi(z_t)
\end{align}
where $\psi_z(z) = \begin{bmatrix}
    z^\top &\varphi_z^\top(z)
\end{bmatrix}^\top$. Moreover, if $r=n$, then the above Koopman operator represents the nonlinear system dynamics up to a diffeomorphism. 
\end{thm}

%% file: Theorems/Proofs/Proof_Theorem3.tex
\begin{pf}
It is given that $dim(\mathcal{O}_y(x))=r$ for $x\in \mathcal{U}\subseteq \mathcal{M}$. Hence, using Theorem \ref{Theorem: Theoretical r-dimensional OC-KO existence}, we can find a Koopman operator of the form (\ref{eq: Koopman operator of Sussman decomposition}) where all the functions of the observation space  $\mathcal{O}_y(x)$ lie in the span of its Koopman observables $\psi_o(x^o)$. From Theorem \ref{Theorem: Theoretical r-dimensional OC-KO existence}, we also know that $x^o$ is formed by a set of $r$ functions in $\mathcal{O}_y(x)$ with a Jacobian of rank $r$ in the neighborhood $\mathcal{U}$. So, we construct the vector \[\begin{bmatrix}
\big{(}h(x)\big{)}^\top & \big{(}h(f(x))\big{)}^\top & \cdots & \big{(}h(f^{n_d-1}(x))\big{)}^\top
\end{bmatrix}^\top\]
and for some $n_d \in \mathbb{N}$, this vector will contain $r$ functions which satisfy the Jacobian criteria for the choice of $x^o$. At time point $n_dt$, this vector becomes the delay embedded output $z_t$ as $y_{n_dt+i} = h(f^{i}(x_{n_dt}))$. Then, all the Koopman observables of (\ref{eq: Koopman operator of Sussman decomposition}) can be written as a function of $z_t$ thereby converting (\ref{eq: Koopman operator of Sussman decomposition}) to a Koopman operator representation of the form $\psi_z(z_{t+1}) = K_z\psi(z_t)$. If $r=n$, then using the discrete-time equivalent of Proposition 3.34 in \cite{nijmeijer1990nonlinear}, we can find a diffeomorphic map between $z_t$ and $x$ in the neighborhood $\mathcal{U}$. Therefore, we can claim that $\psi_z(z_{t+1}) = K_z\psi(z_t)$ captures the full system dynamics up to a diffeomorphism. Hence the proof. 
\hspace*{\fill}~\qed\par
\end{pf}

%% file: Examples/Theoretical_Example_1.tex
\textbf{\subsection{Analytical example to illustrate the theoretical results}}
We consider a nonlinear system with an accurate finite dimensional Koopman operator representation \cite{brunton2016koopman} to illustrate the above theorems. The nonlinear system (\ref{eq: nonlinear system with output})
\begin{align*}
    x_{t+1,1} &= ax_{t,1}\\
    x_{t+1,2} &= bx_{t,2} + \gamma x_{t,1}^2\\
    y_t &= x_{t,2}^2
\end{align*}
has a finite dimensional Koopman operator representation (\ref{eq: Koopman operator system with output})
\begin{align*}
    \begin{bmatrix}
        x_{t+1,1}\\
        x_{t+1,2}\\
        \varphi_1(x_{t+1})\\
        \varphi_2(x_{t+1})\\
        \varphi_3(x_{t+1})\\
        \varphi_4(x_{t+1})
    \end{bmatrix} &= 
    \begin{bmatrix}
        a & 0 & 0 & 0 & 0 & 0\\
        0 & b & \gamma & 0 & 0 & 0\\
        0 & 0 & a^2 & 0 & 0 & 0\\
        0 & 0 & 0 & b^2 & \gamma^2 & 2b\gamma\\
        0 & 0 & 0 & 0 & a^4 & 0\\
        0 & 0 & 0 & 0 & \gamma & b
    \end{bmatrix}
    \begin{bmatrix}
        x_{t,1}\\
        x_{t,2}\\
        \varphi_1(x_{t})\\
        \varphi_2(x_{t})\\
        \varphi_3(x_{t})\\
        \varphi_4(x_{t})
    \end{bmatrix}\\
    y_t &= \varphi_2(x_t)
\end{align*}
where the nonlinear observables are $\varphi_1(x) =x_1^2$, $\varphi_2(x) =x_2^2$, $\varphi_3(x) =x_1^4$ and $\varphi_4(x) =x_1^2x_2$. With $V = \begin{bmatrix}
    0_{3\times 3} & \mathbb{I}_{3\times 3}\\
    \mathbb{I}_{3\times 3} & 0_{3\times 3}
\end{bmatrix}$ and choosing the first 3 observables , we get the minimal Koopman operator representation which captures the output (\ref{eq: Koopman operator of Sussman decomposition}) as
\begingroup
\allowdisplaybreaks
\begin{align*}
    \begin{bmatrix}
        \varphi_2(x_{t+1})\\
        \varphi_3(x_{t+1})\\
        \varphi_4(x_{t+1})
    \end{bmatrix} &= 
    \begin{bmatrix}
        b^2 & \gamma^2 & 2b\gamma\\
        0 & a^4 & 0\\
        0 & \gamma & b
    \end{bmatrix}
    \begin{bmatrix}
        \varphi_2(x_{t})\\
        \varphi_3(x_{t})\\
        \varphi_4(x_{t})
    \end{bmatrix}\\
    y_t &= \varphi_2(x_t).
\end{align*}
\endgroup
Some of the key inferences from the above example are
\begin{itemize}
    \item Choosing $x^o =\begin{bmatrix}
        \varphi_2(x) & \varphi_3(x)
    \end{bmatrix}^\top$, we get the state of the reduced Koopman operator representation (\ref{eq: Koopman operator of Sussman decomposition}) as $\psi_o(x^o) = \begin{bmatrix}
        \varphi_2(x) & \varphi_3(x) & (\varphi_2(x)\varphi_3(x))^{0.5}
    \end{bmatrix}^\top$. Moreover, in the neighborhood $\mathcal{U}$ defined by $x_1,x_2 \in (0,\infty)$, the Jacobian of $x^o$ has a constant dimension of 2.
    \item  The observation space of the nonlinear system comprises
    \begin{align*}
        h(x) &= x_2^2\\
        h(f(x)) &= b^2x_2^2 + 2 b \gamma x_1^2 x_2 + \gamma^2 x_1^4\\
        h(f^2(x)) &= b^4 x_2^2  + (2 a^2 b^2 \gamma + 2 b^3 \gamma) x_1^2 x_2 \\&\quad + (a^4 \gamma^2 + 2 a^2 b \gamma^2 + b^2 \gamma^2) x_1^4 \\
        & \hspace{3mm} \cdots
    \end{align*}
    It can be seen that all the functions of the observation space given by $h(f^{i}(x))$ can be written as a linear combination of $\{h(x),h(f(x)),h(f^2(x))\}$ which has an invertible linear transformation with $\{x_2^2, x_1^4, x_1^2x_2\}$ (if system parameters obey $a^2\neq b$ and $\gamma \neq 2b$). Hence, $\mathcal{O}_y(x)$ lies in the span of $\psi(x)$.
    \item In the same neighborhood $\mathcal{U}$ defined above, we can find another $x^o = \begin{bmatrix}
        \varphi_3(x) & \varphi_4(x)
    \end{bmatrix}^\top$ which leads to $\psi_o(x^o) = \begin{bmatrix}
    \frac{\varphi_4^2(x)}{\varphi_3(x)}    & \varphi_3(x) & \varphi_4(x)
    \end{bmatrix}^\top $ \vspace{3pt}showing that $x^o$ is not unique. 
    \item The Jacobian of the set $\{h(x), h(f(x))\}$ has rank 2 in the neighborhood $\mathcal{U}$ defined by $x_1,x_2 \in (0,\infty)$. Hence, we can define a delay embedded output coordinate $z_t = \begin{bmatrix} y_{2t} & y_{2t+1}\end{bmatrix}^\top$
    which a Koopman operator of the form 
    \begin{align*}
        \psi(z_{t+1}) &= K_z\psi(z_t)\\
        \text{where }\psi(z_t) &=\begin{bmatrix} y_{2t} & y_{2t+1} & \varphi_z(y_{2t},y_{2t+1})
    \end{bmatrix}^\top.
    \end{align*}
    The computation is straight-forward and lengthy. Hence, we just state the result in abstraction. The delay embedded output can capture the state dynamics up to a diffeomorphism.
\end{itemize}
This analytical example illustrates all of the above theoretical results.

%% file: sections/5_SimulatedResults.tex
\section{Simulation Results}\label{sec: Simulation Results}
In this section, we demonstrate that the theory in Section \ref{sec: Main Results} can be used in complex nonlinear systems to determine the critical states that drive an output performance objective. Specifically, in biological systems, we tackle an important problem --- what are the genes (state) that affect a certain phenotype (output performance metric)? 

For each system, we start by learning Koopman operator representations with output equations (\ref{eq: Koopman operator system with output}) using OC-deepDMD algorithm as mentioned in Appendix \ref{appendix: ocdeepDMD} to capture the nonlinear dynamics of the form (\ref{eq: nonlinear system with output}). For each learned model, we compute its 1-step and n-step prediction accuracy of both states and outputs to ensure that the model captures the nonlinear dynamics with a high accuracy.
For the 1-step prediction, given the state at one time point, we predict the next time point. For the n-step prediction, given only the initial condition of the state, we predict the states for all future time points within the time period of the simulation run. The 1-step prediction accuracy is a representation of how well we solve the OC-deepDMD optimization problem (Section \ref{sec: OC-DMD}) and the n-step prediction accuracy is a representation of how well we capture the actual nonlinear dynamics of the system. To compute the accuracy of the predictions with the true data, we use the $r^2-$score, also called the coefficient of determination. 

Once we have learned a linear Koopman model (\ref{eq: Koopman operator system with output}) that adequately captures the system dynamics, we reduce this model to the minimal Koopman model (\ref{eq: Koopman operator of Sussman decomposition}) that drives the output performance measure. This procedure is highlighted in Appendix \ref{appendix: observable decomposition}. The Koopman observables $\psi_o(x)$ of the model (\ref{eq: Koopman operator of Sussman decomposition}) capture the full observation space $\mathcal{O}_y(x)$ but the more important information is how much do each of the state variables in $x$ contribute to $\psi_o(x)$. We have developed a sensitivity computation algorithm as a followup to the OC-deepDMD algorithm to identify the genes in order of their importance to the given output performance measure.

\input{Images/FigureTexFiles/Example1}

\subsection{Example 1 - Finding Critical Genes to Control Bacteria Growth}
One of the prominent performance metrics (phenotype) used in biological systems is the population growth of cell cultures. Specifically, the challenge is to identify genes that are responsible for the cells to proliferate when subject to different growth substrates like sugars, proteins and other conditions like pH and oxygen levels. We illustrate this challenge by simulating the ccd antitoxin-toxin system \cite{aghera2020mechanism} which is known to regulate growth in bacteria. Specifically, the dynamic interaction of CcdA antitoxin and CcdB toxin regulates the concentration of DNA gyrase which plays a crucial role in relieving the topological stress while the DNA is transcribed by the RNA-polymerase enzyme. DNA gyrase complex enhances the production of proteins from cellular DNA which could either up-regulate or down-regulate the cell proliferation process. We simulate a simplified model of the complex network using the ccd antotoxin-toxin reaction network in \cite{aghera2020mechanism} with its output DNA gyrase modulating the expression of four genes, which directly impact (positively or negatively) the growth output following Monod's growth kinetics model \cite{luong1987generalization}. The gene network is shown in Fig. \ref{fig: example1}-A and the system dynamics along with the simulation details are given in Appendix \ref{sec: Simulation example 1}.

Given that we have the state and the output data of the above nonlinear system, our objective to identify the states (genes) that impact the output (growth performance metric) dynamics. We begin by learning the Koopman operator with output (\ref{eq: Koopman operator system with output}) using the OC-deepDMD algorithm (Appendix \ref{appendix: ocdeepDMD}) and the model prediction on a random initial condition is shown in Figure \ref{fig: example1}(b). For a test data set, the identified optimal model has a state ($x$) prediction accuracy of  $98.6\%$ for 1-step and $98.4\%$  for n-step and an output ($y$) prediction accuracy of $99\%$ for 1-step and $82\%$ for n-step. On using the linear observable decomposition procedure from Appendix \ref{appendix: observable decomposition}, we can reduce the identified 24-dimensional Koopman operator model (\ref{eq: Koopman operator system with output}) to a 15-dimensional minimal Koopman model (\ref{eq: Koopman operator of Sussman decomposition}) with Koopman observables $\psi_o(x)$ that capture the output. We evaluate the sensitivities of each of the functions in $\psi_o(x)$ with respect to the base states $x$ as described in Appendix \ref{appendix: sensitivity}; the sensitivity matrix is shown in lower heatmap plot of Fig. \ref{fig: example1}(c) and the contribution of each state to $\psi_o(x)$ (the Euclidean norm of the sensitivity matrix) is shown in the upper bar plot of Fig. \ref{fig: example1}(c). 

The following results reveal the success of our algorithm:
\begin{enumerate}
    \item States $x_8$ and $x_{11}$(red), that directly impact the output, have the most contribution towards $\psi_o(x)$.
    \item States $x_9$ and $x_{10}$(brown) which have no impact on the output provide the least contribution to $\psi_o(x)$.
    \item States $x_1$ through $x_7$(blue and green) which represent the toxin-antitoxin system are the secondary states that indirectly contribute to the output. It can be seen that their contributions to $\psi_o(x)$ lie between the two extreme cases and their contribution to $\psi_o(x)$ reduces as the gene is located farther away in the network from the output (as we see in transitioning from the genes in green to the genes in blue).
\end{enumerate}
It is evident that the results are not perfect (like $x_8$ and $x_{11}$ have non-zero contributions to $\psi_o(x)$). The imperfections are a result of various numerical approximations. The state-inclusive Koopman model is a finite dimensional approximation learned by minimizing the 1-step prediction error and naturally, by its very formulation, it cannot capture the entirety of the nonlinear dynamics. The linear observable decomposition in Appendix \ref{appendix: observable decomposition} is numerically approximated.  Despite the various sources of error, to a large extent, our algorithm can order the genes (states) based on their importance to the output performance measure.

\begin{figure*}[htp!]
\centering
    \includegraphics[width=\textwidth,keepaspectratio]{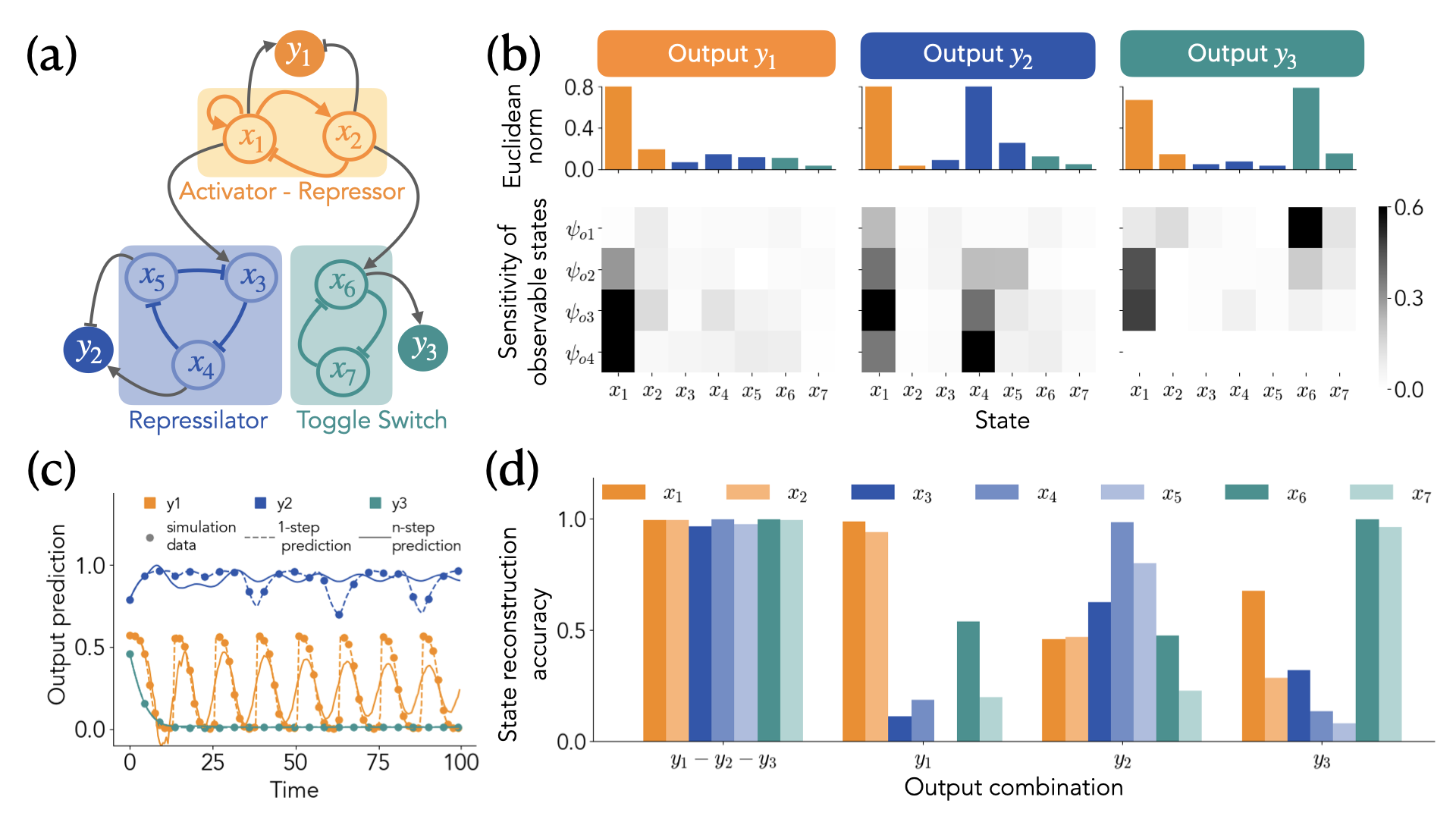}
    \caption{
    \textbf{Example 2 - Finding Critical Genes In Composed Genetic Circuit Networks:}  \textbf{(a)} A complex genetic circuit formed by interconnecting three well-studied genetic circuits with an output measured from each of the core circuits \textbf{(b)} The sensitivity heat map (lower) of the Koopman observables $\psi_o(x)$ of the minimal Koopman operator (\ref{eq: Koopman operator of Sussman decomposition}) for each output and the corresponding 2-norm bar plot (upper) showing the contributions of each gene to the Koopman observables $\psi_o(x)$
    \textbf{(c)} The 1-step and n-step prediction of the optimal delay embedded Koopman operator (\ref{eq: delay embedded Koopman operator}) learned using only the output data $y_1$, $y_2$, and $y_3$
    \textbf{(d)} State reconstruction accuracy from the Koopman observables of optimal delay embedded Koopman operators for various combinations of outputs.
    }\label{fig: example2}
\end{figure*}

\textbf{\subsection{Example 2 - Finding Critical Genes In Composed Genetic Circuit Networks}}
 We consider another complex genetic circuit composed of three interconnected subsystems (taken from \cite{del2015biomolecular}): the activator-repressor, the repressilator, and the toggle switch with a single output measured from each subsystem as shown in Fig. \ref{fig: example2}(a). The nonlinear system dynamics and simulation details are given in Appendix \ref{sec: Simulation System 2}. We show that the observability of Koopman operators can reveal the genes that impact each individual output. 

\subsubsection{Trade-off in learning state-inclusive Koopman operator models\textbf{:}}
We adopt the same methodology as in Example 1. We begin by learning a Koopman operator model with output (\ref{eq: Koopman operator system with output}) as mentioned in Appendix \ref{appendix: ocdeepDMD}. The model has a state $x$ prediction accuracy of $99.8\%$ for 1-step and $69\%$ for n-step predictions and an output $y$ prediction accuracy of $98\%$ for 1-step and $64\%$ for n-step predictions. The low n-step prediction accuracy is a consequence of the trade-off  between the  n-step prediction accuracy and the state-inclusivity of the Koopman operator representation; the state-inclusive Koopman operator enables easiest reconstruction of the original state $x$ of the nonlinear system (by simply dropping the nonlinear Koopman observables) but the state-inclusive Koopman operator model converges to a single equilibrium point (by construction) which is not suitable for systems that exhibit an oscillatory steady state response like the nonlinear system under consideration. Moreover, the OC-deepDMD objective function is constructed to minimize only the 1-step predictions and hence, does not guarantee n-step prediction accuracy. We see that the Koopman model learned from OC-deepDMD algorithm that minimizes only 1-step prediction error is still adequate for gene identification.

\subsubsection{Linear observable decomposition of the Koopman model reveals the critical genes that impact each output\textbf{:}}
For each output in the vector of outputs, we consider the row of $W_h$ corresponding to that output and learn the minimal Koopman operator (\ref{eq: Koopman operator of Sussman decomposition}) that captures the dynamics of that output using the method in Appendix \ref{appendix: observable decomposition} and identify the sensitivity matrices and the Euclidean norm as in Example 1 using the approach in Appendix \ref{appendix: sensitivity}. The sensitivity plots are shown in Fig. \ref{fig: example2}(b).

From the genetic circuit Fig \ref{fig: example2}(a) and the data based ordering of state contributions in Fig \ref{fig: example2}(b), we can see that the key results are captured
\begin{enumerate}
    \item The output $y_1$ is mainly influenced by $x_1$ which activates $y_1$ followed by $x_2$ which represses $y_1$
    \item The output $y_2$ is mainly influenced by $x_4$ which activates $y_2$, followed by $x_1$ which activates the repressilator and $x_5$ which represses $y_2$.
    \item The output $y_3$ is mainly impacted by $x_6$ and $x_1$ followed by $x_7$ and $x_2$.
\end{enumerate}
In addition to the main results, there are residual contributions by each state to each output. This can be attributed to three sources of error: the low n-step prediction accuracy, the numerical approximations and linear correlations between the state variables.  

An important observation across Examples 1 and 2 is that both activator ($x_8$) and repressor ($x_{11}$) genes in Example 1 are recognized as significant genes whereas in Example 2, the significance of activator genes ($x_1,x_4,x_6$) is more prominent than the repressor genes ($x_2,x_5$). The explanation for the same lies in how much the activators and repressors impact the outputs of the system. In Example 1, we can see that that activator $x_8$ and repressor $x_{11}$ both directly impact the output and both growth and decaying effects are captured in the output. In Example 2, the output $y_3$ has no direct repressors impacting it and though the outputs $y_1$ and $y_2$ are repressed by genes $x_2$ and $x_5$ respectively, the effect of repression is not prominent as witnessed by the absence of any decaying effects in the evolution of the outputs $y_1$ and $y_2$. Hence, it is evident that the algorithm captures the important genes based on how much influence the genes have on the phenotype and not just the proximity of the genes to the phenotype (in the gene network). 

\subsubsection{The fusion of the three outputs contain adequate information to represent the full state of the system\textbf{:}} We consider the possibility that there might not be adequate information in the output measurements to inform the genes. To ensure that the outputs are rich enough to capture the state information, we make use of Theorem \ref{Theorem: delay embedded Koopman operator}; we identify a delay embedded Koopman operator model using only the output measurements and examine if a diffeomorphic map exists between the observables of the delay embedded Koopman operator model and the state $x$. The delay embedded output is given by 
\[
z_{t} = \begin{bmatrix} y_{n_dt}^\top & y_{(n_d+1)t}^\top & \cdots y_{(2n_d-1)t}^\top \end{bmatrix}^\top
\]
and the delay embedded Koopman operator is solved by adopting the same method as in Appendix \ref{appendix: ocdeepDMD} except with a different objective function
\[
\min_{\psi,K} ||\psi(Z_F^{train}) - K\psi(Z_P^{train})||_F^2
\]
which has an added hyperparameter, $n_d$. The parameter $n_d$ is the number of output delay embeddings used to construct the observables $\psi(z)$ of the delay embedded Koopman operator $K$. 

The predictions of the delay-embedded Koopman operator model (with optimal $n_d=4$) on a random initial condition (from test data set) is shown in Figure \ref{fig: example2}(c). The optimal delay embedded Koopman model using all the outputs ($y_1,y_2$ and $y_3$) has a $99.8\%$ 1-step prediction accuracy and $58.1\%$ n-step prediction accuracy. As an outcome of Theorem \ref{Theorem: delay embedded Koopman operator}, we know that the Koopman observables $\psi(z)$ capture the entire observation space. Hence, if the output measurements capture the full system dynamics, we should be able to find a diffeomorphic map between $\psi(z)$ and $x$. We learn a numerical diffeomorphic map using the method in Appendix \ref{appendix: diffeomorphic map}. Then, we use the numerical diffeomorphic map to reconstruct the state and the reconstruction accuracy of each state is shown as a bar plot in Figure \ref{fig: example2}(d) above $y_1-y_2-y_3$. We see that by using all the outputs, all the states can be almost accurately reconstructed.  When we repeat the same process using single measurements like $y_1$, $y_2$ and $y_3$, we see that only partial states show accurate reconstruction. Therefore, we conclude that there is adequate information in the output measurements to capture each gene and our sensitivity analysis orders the genes by how much impact they have on a specified output.

Through the simulation examples, we see that the observability of linear high-dimensional Koopman operator models with linear output equations can be used as a proxy to the observability of the nonlinear systems. In biological systems, we see that this approach is very useful to discover genes that drive various phenotypic behaviors.

%% file: Images/FigureTexFiles/Example1.tex
\begin{figure*}[htp!]
\centering
    \includegraphics[width=\textwidth,keepaspectratio]{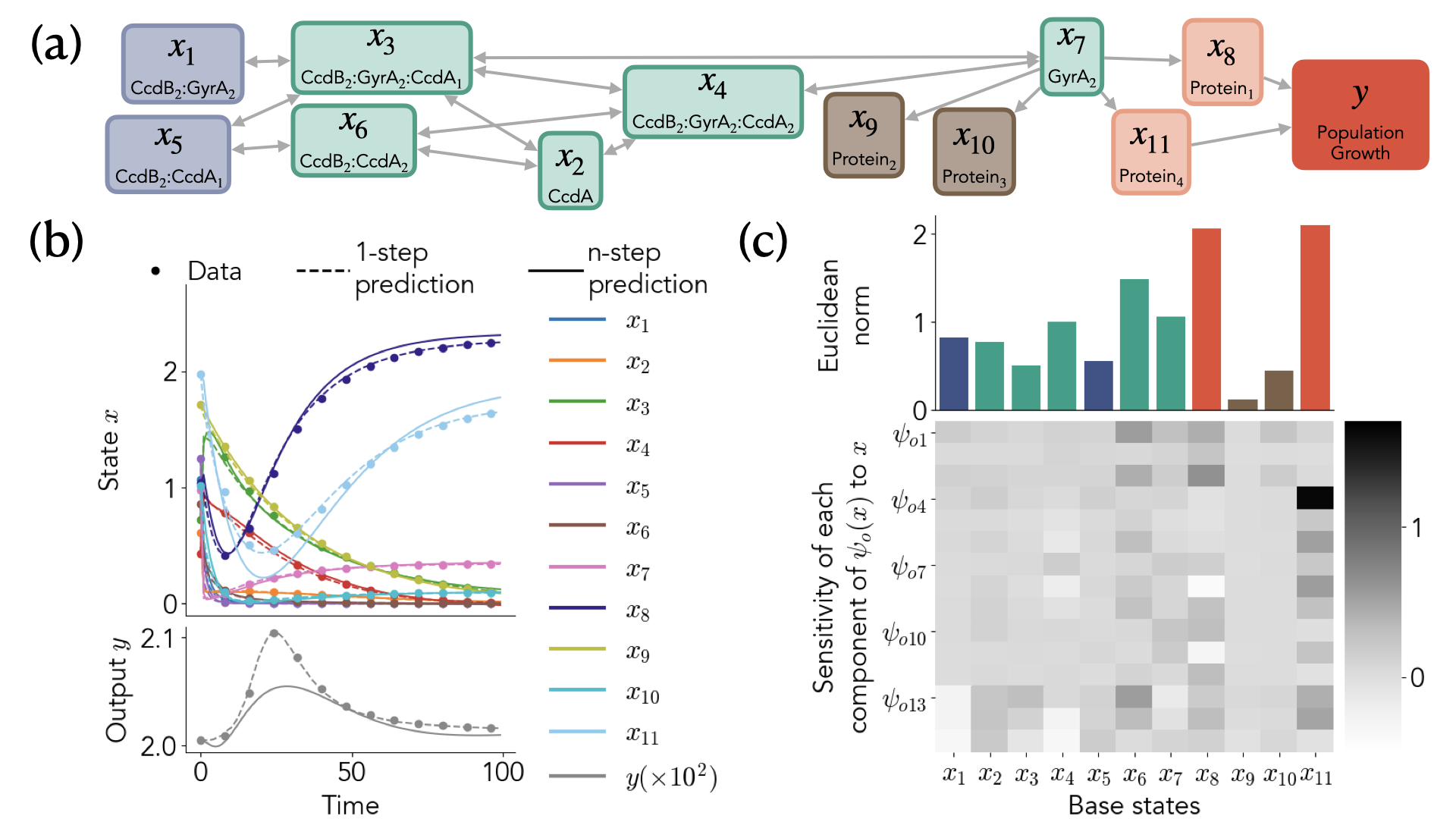}
    \caption{\textbf{Example 1 - Finding Critical Genes to Control Bacteria Growth:} \textbf{(a)} The directed graph of the reaction network: : the states x1 through x7 ($x_1$ and $x_5$ indicated in blue initiate the network and  $x_2$, $x_3$, $x_4$, $x_6$, and $x_7$ indicated in green are the intermediates and products of the reaction network) is the toxin-antitoxin system taken from  \cite{aghera2020mechanism} , states $x_8$ and $x_{11}$ (light red) are two proteins enhanced by the gyrase enzyme that have a positive and negative effect on growth output (dark red) respectively, and states $x_9$ and $x_{10}$(brown) are two proteins enhanced by gyrase enzyme but have no association with growth \textbf{(b)} The 1-step and n-step predictions of the Koopman operator model with output (\ref{eq: Koopman operator system with output}) learned using the OC-deepDMD algorithm in Appendix \ref{appendix: ocdeepDMD} \textbf{(c)} The bottom heat map is the sensitivity of the Koopman observables $\psi_o(x)$ of the minimal Koopman operator (\ref{eq: Koopman operator of Sussman decomposition}) with respect to the system states $x$. The top bar plot computes the euclidean norm for each column in the sensitivity matrix to represent the relative contributions of each state $x$ to $\psi_o(x)$. 
    }\label{fig: example1}
\end{figure*}

%% file: sections/6_Conclusion.tex
\section{Conclusion}\label{sec: Conclusion}
In this paper, we show how linear observability analysis of Koopman operator models ties to the observability analysis of nonlinear dynamical systems. We provide algorithms to learn Koopman operators which constrain the outputs to lie in the span of Koopman observables. We show how decomposition of these output-inclusive Koopman operator models discover a reduced set of Koopman observables that drive output dynamics. The techniques can be seamlessly applied to other complex systems involving data-driven learning of governing dynamics.

In biological systems, we show how to find the minimal Koopman operator models that capture the output dynamics and use sensitivity analysis to discover the genes (states) that drive a phenotypic behavior (output performance metric). Through this work, we solve the first step toward our ultimate objective of controlling the expression of critical genes that regulate phenotypic behavior in biological systems.

%% file: sections/7_Acknowledgement.tex
\begin{ack}   
The authors would like to graciously thank Igor Mezic, Nathan Kutz, Joshua Elmore, Adam Deutschbauer, Bassam Bamieh, and Charles Johnson for insightful discussions. Any opinions, findings, conclusions, or recommendations expressed in this material are those of the authors and do not necessarily reflect the views of the Defense Advanced Research Project Agency, the Department of Defense, or the United States government.    

This work was supported in part by the Department of Energy’s Biological and Environmental Research office, as a contribution of the Pacific Northwest National Laboratory Science Focus Area ``Persistence Control of Engineered Functions in Complex Soil Microbiomes'' via PNNL subcontract numbers 545157 and 490521.

This work  was also partially supported by DARPA, AFRL under contract numbers FA8750-17-C-0229, HR001117C0092, HR001117C0094, DEAC0576RL01830, as well as funding from the Army Research Office's Young Investigator Program under grant number W911NF-20-1-0165.   Supplies for this work were partially supported by the Institute of Collaborative Biotechnologies, via grant W911NF-19-D-0001-0006.

\end{ack}

%% file: sections/8_Appendix.tex
\section{Appendix}\label{sec: appendix}

\subsection{Learning a Koopman operator model for the nonlinear dynamical system with outputs: \textit{output-constrained deep dynamic mode decomposition (OC-deepDMD) algorithm }} \label{appendix: ocdeepDMD}

\subsubsection{Data generation} For a given nonlinear system, we simulate the nonlinear model for multiple initial conditions and record both the states $x$ and the outputs $y$. The data is equally split among training, validation and test sets. To ensure equal representation of data across the three sets, the initial conditions are randomly sampled from a uniformly distributed bounded phase space. 
\subsubsection{Data preprocessing} For each initial condition ($i$) in the training, validation and test datasets, the generated data is organized as
\begin{align*}
    X_p^{(i)} &= \begin{bmatrix}x^{(i)}_{0} & x^{(i)}_{1} & \cdots & x^{(i)}_{N_{sim}-1}  \end{bmatrix}\\
    X_f^{(i)} &= \begin{bmatrix}x^{(i)}_{1} & x^{(i)}_{2} & \cdots & x^{(i)}_{N_{sim}}  \end{bmatrix} \\
    \text{ and } Y_p^{(i)} &= \begin{bmatrix}y^{(i)}_{0} & y^{(i)}_{1} & \cdots & y^{(i)}_{N_{sim}-1}  \end{bmatrix}
\end{align*}
where $x|y_t^{(i)}$ indicates either the state $x$ or the output $y$ at time point $t$ generated from the $i^{th}$ initial condition. The data across the snapshots are concatenated together as $S= \begin{bmatrix}S^{(1)} & S^{(2)} & \cdots \end{bmatrix}$ where $S=X_p, X_f$ or $Y_p$. The training data $X_p^{train}$ and $Y_p^{train}$ are used to identify the mean and standard deviation of each variable and all the data are standardized (subtracted by the computed mean and divided by the computed standard deviation). 
\subsubsection{Learning an optimal model for a set of hyperparameters} We use tensorflow in Python to setup neural networks, the outputs of which represent the observables of the Koopman operator that we want to learn. The hyperparameters of the model include the number of nodes in each layer of the neural network, the number of layers in the neural network, the activation function in each node and the number of output nonlinear observables ($\varphi(x)$). We append the nonlinear observables $\varphi(x)$ to the states $x$ and the bias term $1$ to avoid trivial solutions and the full observable vector at a single time point is given by $\psi(x) = \begin{bmatrix}x^\top & \varphi^\top(x) & 1\end{bmatrix}^\top$. We also initialize the matrices $K$ and $W_h$ from (\ref{eq: Koopman operator system with output}) in the tensorflow environment, set up the objective function 
\begin{align*}
    \min_{\psi,K,W_h} ||\psi(X_F^{train}) &- K\psi(X_P^{train})||_F^2 \\&+ ||Y_P^{train}-W_h\psi(X_P^{train})||_F^2
\end{align*}
and use Adagrad optimizer in Python to implement stochastic gradient descent with various step sizes to identify an optimal model for a given set of hyperparameters. 

\subsubsection{Learning model with optimal hyperparameters} For various combination of the hyperparametes, we learn an optimal Koopman operator model. We evaluate 1-step and n-step state and output prediction accuracy for each model across the training and validation datasets:
\begin{align*}
    r^2_{s,(1|n)-step} &= 1 - \frac
    {\sum_{i}\sum_{j} (s_j^{(i)}-\hat{s}_j^{(i)})^\top(s_j^{(i)}-\hat{s}_j^{(i)})}
    {\sum_{i}\sum_{j} (s_j^{(i)}-\bar{s}_j^{(i)})^\top(s_j^{(i)}-\bar{s}_j^{(i)})}
\end{align*}
where $s$ is either the state $x$ or the output $y$, $j$ indicates the time point and $i$ indicates the initial condition the data is generated from. $\bar{s}$ is the mean of $X_p^{train}$ for state $x$ and mean of $Y_p^{train}$ for output $y$ and $\hat{s}_j^{(i)}$ is the inverse standardization of 
\begin{itemize}
    \item $\begin{bmatrix}\mathbb{I}_n & 0\end{bmatrix}K\psi(x_{j-1}^{(i)})$ for 1-step $x$ prediction,
    \item $W_hK\psi(x_{j-1}^{(i)})$ for 1-step $y$ prediction,
    \item $\begin{bmatrix}\mathbb{I}_n & 0\end{bmatrix}K^j\psi(x_{0}^{(i)})$ for n-step $x$ prediction, and
    \item $W_hK^j\psi(x_{0}^{(i)})$ for n-step $y$ prediction.
\end{itemize}
We use these metrics to settle on a model that is optimized in both parameters and hyperparameters.

\subsection{Learning the observable decomposition form of a Koopman operator model with output}\label{appendix: observable decomposition}
Given that we have a state-inclusive Koopman operator model of the form (\ref{eq: Koopman operator system with output}) identified using the method in Appendix  \ref{appendix: ocdeepDMD}, we want to find a dimensionality reduced model of the form (\ref{eq: Koopman operator of Sussman decomposition})--- a model with minimal Koopman observable functions to capture the output dynamics. In practise, (\ref{eq: Koopman operator system with output}) is typically a finite dimensional approximation. We find the observability matrix of the identified Koopman system $\mathcal{O}_y(x) = \begin{bmatrix}W_h^\top & (W_hK)^\top & \cdots &  (W_hK^{n_L})^\top \end{bmatrix}^\top$ and its right singular vectors ($V$). Then we can transform (\ref{eq: Koopman operator system with output}) as $\psi_{ou}(x) = V^\top\psi(x)$, $\tilde{K} = V^\top K V$ and $\tilde{W}_h = W_h V$. In theory, the upper right block of $\tilde{K}$  and the right block of $\tilde{W}_h$ should be 0 (as seen in Corollary \ref{cor: linear observable decomposition}). Due to numerical approximation, perfect zero cannot be obtained. The challenge is to estimate the dimension of $\psi_o(x)$ ($n_{oL}$) in (\ref{eq: Koopman operator of Sussman decomposition}) where  $\psi_o(x)$ is the first $n_{oL}$ elements of $\psi_{ou}(x)$. We use the property that $\psi_o(x)$ can accurately capture the output; we increase $n_{oL}$ from $1$ to $n_L$ and examine at what value of $n_{oL}$ can $\psi_o(x)$ capture 99\% ($r^2$ score) of the output predicted by (\ref{eq: Koopman operator system with output}). This yields the required reduced model of the form (\ref{eq: Koopman operator of Sussman decomposition}) with the required properties intact.

\subsection{Computing the sensitivity of each nonlinear function in $\psi_o(x)$ with respect to the base coordinate states $x$} \label{appendix: sensitivity}
Neural networks are typically used to approximate functions. In the minimal Koopman operator that captures the output dynamics, the set of nonlinear observable functions $\psi_o(x)$ is captured by a neural network that we implement using tensorflow in python. To compute the sensitivity of a single function in the set $\psi_o(x)$ with respect to a single state variable in $x$, we simply use the \textit{gradients} function in tensorflow package of python. We evaluate the gradient at all training data points and store the maximum. We evaluate this maximum sensitivity for each function in $\psi_o(x)$ with respect to each state variable in $x$. This yields a matrix which is not very intuitive to interpret which state is more important. So, we compute the euclidean norm of the sensitivity matrix for each state variable in $x$ across the maximum sensitivities of all functions in $\psi_o(x)$ with respect to that state variable in $x$. 

\subsection{Learning the differomorphic map between the delay embedded output $z$ and the base coordinate state $x$}\label{appendix: diffeomorphic map}
\begin{figure}[H]
\centering
    \includegraphics[width=3in,keepaspectratio]{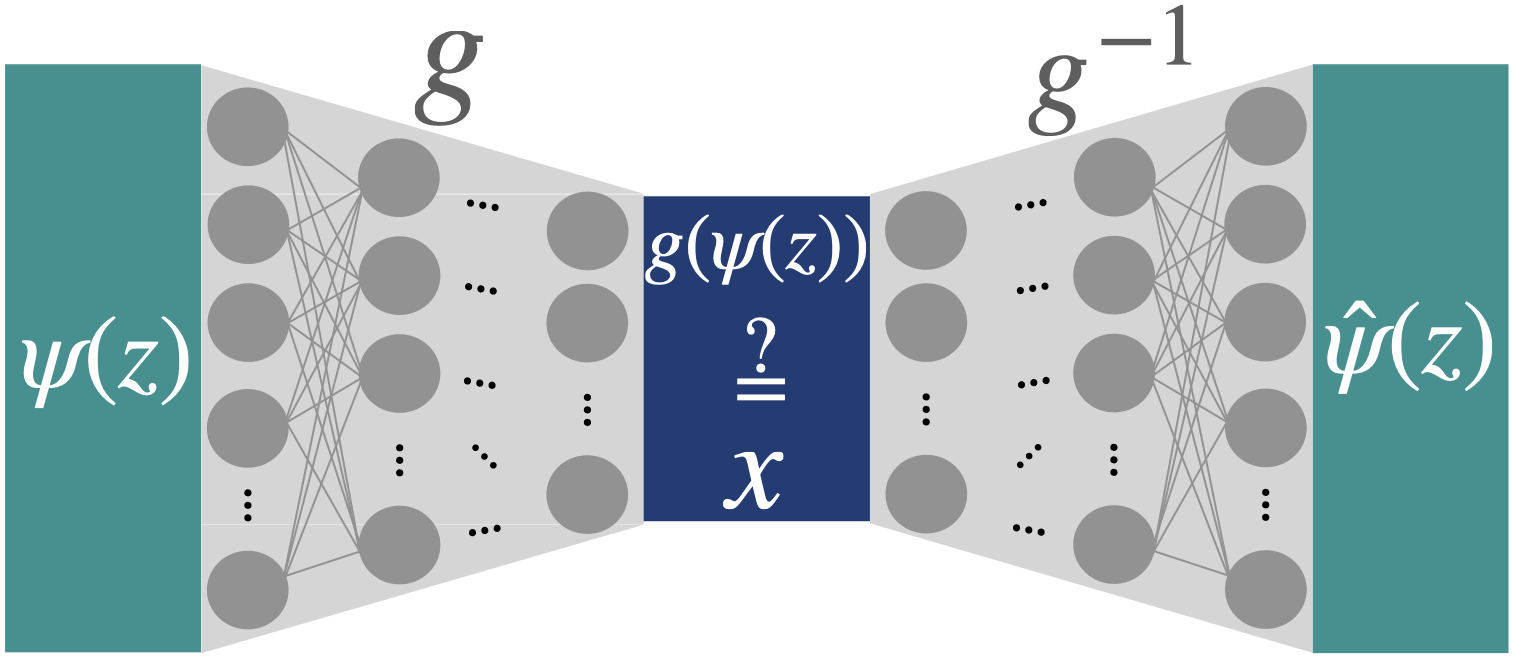}
\end{figure}

Given the delay embedded output $\psi(z)$ and the state $x$, we represent the diffeomorphic map $\big{(}$the forward transform $g:\psi(x) \rightarrow x$ and the inverse transform $g^{-1}:x \rightarrow \psi(x)\big{)}$ using the autoencoder-decoder neural network. Specifically, we formulate the multi-objective optimization problem
\begin{align*}
    \min_{g,g^{-1}} ||\psi(z)  - g^{-1}(g(\psi(z)))||_F^2 + ||x  - g(\psi(z))||_F^2
\end{align*}
in Python using Tensorflow and solve it by using the Adagrad optimizer to implement stochastic gradient descent. The two objectives that the above optimization targets are \textit{(i)} to transform $\psi(z)$ to a reduced coordinate space and \textit{(ii)} to get the reduced coordinates close to the state $x$ as much as possible. 

\subsection{Simulation parameters for Example 1}\label{sec: Simulation example 1}
The dynamics of the gene network with growth output  in Fig. \ref{fig: example1}(a) is given by:

\begin{align*}
    \dot{x}_1 &= -k_{1f}x_1x_2 +k_{1r}x_3 - \gamma_1x_1 + u_0\\
    \dot{x}_2 &= -k_{1f}x_1x_2 +k_{1r}x_3 -k_{2f}x_2x_3 + k_{2r}x_4 - k_{5f}x_2x_5\\
    &\quad + k_{5r}x_6 - \gamma_2x_2\\
    \dot{x}_3 &= k_{1f}x_1x_2 - k_{1r}x_3 -k_{2f}x_2x_3 + k_{2r}x_4 - k_{4f}x_3\\
    &\quad + k_{4r}x_5x_7 - \gamma_3x_3\\
    \dot{x}_4 &= k_{2f}x_2x_3 - k_{2r}x_4 - k_{3f}x_4 + k_{3r}x_6x_7 - \gamma_4x_4\\
    \dot{x}_5 &= k_{4f}x_3 - k_{4r}x_5x_7 - k_{5f}x_2x_5 + k_{5r}x_6- \gamma_5x_5\\
    \dot{x}_6 &= k_{5f}x_2x_5 - k_{5r}x_6 + k_{3f}x_4 - k_{3r}x_6x_7- \gamma_6x_6\\
    \dot{x}_7 &= k_{3f}x_4 - k_{3r}x_6x_7 + k_{4f}x_3 - k_{4r}x_5x_7- \gamma_7x_7\\
    \dot{x}_8 &=\frac{a_1(x_7/k_1)^{n_1}}{1+(x_7/k_1)^{n_1}} - d_1x_8\\
    \dot{x}_9 &=\frac{a_2(x_7/k_2)^{n_2}}{1+(x_7/k_2)^{n_2}} - d_2x_9\\
    \dot{x}_{10} &=\frac{a_3(x_7/k_3)^{n_3}}{1+(x_7/k_3)^{n_3}} - d_3x_{10}\\
    \dot{x}_{11} &=\frac{a_4(x_7/k_4)^{n_4}}{1+(x_7/k_4)^{n_4}} - d_4x_{11}\\
    y &= y_o exp\Big{(}\frac{\mu_yx_8}{K_y + x_8 + x_{11}}\Big{)}.
\end{align*}

The simulation parameters of the system are
$k_{1f} = 1.4 M^{-1}s^{-1}$, $k_{1r} =0.003s^{-1}$, $k_{2f} =1.1M^{-1}s^{-1}$, $k_{2r} =0.19s^{-1}$, $k_{3f} =0.04s^{-1}$, $k_{3r} =2.2M^{-1}s^{-1}$, $k_{4f} =0.0035s^{-1}$, $k_{4r} =2.2M^{-1}s^{-1}$, $k_{5f} =0.14M^{-1}s^{-1}$, $k_{5r} =0.13s^{-1}$,  $a_1=0.8Ms^{-1}$, $k_1=0.3M$, $n_1=2$, $a_2=1.9Ms^{-1}$, $k_2=2M$, $n_2=5$, $a_3=4Ms^{-1}$, $k_3=4M$, $n_3=2$, $a_4=0.7Ms^{-1}$, $k_4=0.5M$, $n_4=3$, $\gamma_1 =0.3s^{-1}$, $\gamma_2 =0.1s^{-1}$, $\gamma_3 =0.03s^{-1}$, $\gamma_4 =0.02s^{-1}$, $\gamma_5 =0.4s^{-1}$, $\gamma_6 =0.09s^{-1}$, $\gamma_7 =0.01s^{-1}$, $d_1 =0.2s^{-1}$, $d_2 =0.03s^{-1}$,  $d_3 =0.3s^{-1}$, $d_4 =0.1s^{-1}$, $\mu_y =10$, $y_0=0.02$, $K_y=10M$. The initial condition of the state is $x_0 = \begin{bmatrix} 0.4, 0.1, 0.2, 0.4, 0.3, 0.8, 0.5, 0.3, 0.8, 0.1, 1.8\end{bmatrix}^\top + e$ where $e\in \mathbb{R}^{11 \times 1}$ with each entry in $e$ uniformly distributed in the range $[0,1]$. The simulation time ($T_s$) is $1s$ with a simulation of $100s$ for each initial condition.

\subsection{Simulation System 2}\label{sec: Simulation System 2}
The gene network dynamics in Fig. \ref{fig: example2}(a) is modeled as:
\begin{subequations}\label{eq: AR_RRR_TS}
    \begin{align}
        \intertext{Activator repressor}
        \dot{x}_1 &= \frac{\kappa_{1}}{\delta_{1}}. 
        \frac{\alpha_1(x_1/K_{1})^{n_{1}}+ \beta_{1}}{1+ (x_1/K_{1})^{n_{1}}+(x_2/K_{2})^{m_{1}}} - \gamma_{1}x_1\nonumber\\
        \dot{x}_2 &= \frac{\kappa_{2}}{\delta_{2}}. 
        \frac{\alpha_2(x_1/K_{1})^{n_{1}}+ \beta_{2}}{1+ (x_1/K_{1})^{n_{1}}} - \gamma_{2}x_2\label{eq: AR}\\
        y_1 &= V_1 \frac{(x_1/K_{11})^{n_{11}}}{1 + (x_1/K_{11})^{n_{11}} + (x_2/K_{12})^{n_{12}}}\nonumber
        \intertext{Repressilator}
        \dot{x}_3 &= c_{1,3}x_1 +\frac{\alpha_{3}}{1 + (x_5/K_{5})^{n_{5}}} - \gamma_{3}x_3\nonumber\\
        \dot{x}_4 &= \frac{\alpha_{4}}{1 + (x_3/K_{3})^{n_{3}}} - \gamma_{4}x_4\label{eq: RRR}\\
        \dot{x}_5 &= \frac{\alpha_{5}}{1 + (x_4/K_{4})^{n_{4}}} - \gamma_{5}x_5\nonumber\\
        y_2 &= V_2 \frac{(x_4/K_{24})^{n_{24}}}{1 + (x_4/K_{24})^{n_{24}} + (x_5/K_{25})^{n_{25}}}\nonumber\\
        \intertext{Toggle switch}
        \dot{x}_6 &= c_{2,6}x_2 + \frac{\alpha_{6}}{1 + (x_7/K_{6})^{n_{6}}} - \gamma_{6}x_6\nonumber\\
        \dot{x}_7 &= \frac{\alpha_{6}}{1 + (x_6/K_{6})^{n_{6}}}- \gamma_{7}x_7 \label{eq: TS}\\
        y_3 &= V_3 \frac{(x_6/K_{36})^{n_{36}}}{1 + (x_6/K_{36})^{n_{36}}}\nonumber
    \end{align}
\end{subequations}
The parameters of the activator reperssor 
$\kappa_{1} = 1$, $\delta_{1} =1$, $\alpha_1 = 250$, $K_1 =1$, $n_1 =2$, $\beta_1 =0.04$, $K_2 = 1.5$, $m_1 =3$, $\gamma_1 =1$, $\kappa_{2} = 1$, $\delta_{2} =1$,  $\alpha_2 = 30$,  $\beta_2 = 0.004$, $\gamma_2 =0.5 $, $V_1 =2$, $K_{11}=1$, $n_{11}=1$, $K_{12} =0.4$ and $n_{12}=1$. The parameters of the repressilator are $c_{1,3}=0.1$, $\alpha_3 =$, $\alpha_4 =$, $\alpha_5 =$, $\gamma_3 =0.3$, $\gamma_4 =0.3$, $\gamma_5 =0.3$, $K_5=1$, $n_5=2$, $K_3=1$, $n_3=4$, $K_4=1$, $n_4=3$, $K_{24}=0.02$, $n_{24}=1$, $K_{25}=1$, $n_{25}=2$ and $V_2=1$. The parameters of the toggle switch are $c_{2,6}=0.001$, $\alpha_6=1$, $K_6=10$, $n_6=1$, $\gamma_6=0.09$, $\gamma_7=0.09$, $K_{36}=120$, $n_{36}=1$ and $V_3=1$. The sampling time was $0.5s$ and the duration of each simulation was $100s$. The initial conditions for the simulation are $x_0 =[100.1,20.1,10.,10.,10.,100.1,100.1]^\top + e$ where $e\in \mathbb{R}^{7 \times 1}$ uniformly distributed in $[0,4]$.